%% file: agt-4-26.tex
\documentclass{gtart_h}

\input agtout

\lognumber{26}
\volumenumber{4}
\volumeyear{2004}
\papernumber{26}
\published{8 August 2004}
\pagenumbers{571}{594}
\received{2 January 2003}
\revised{1 June 2004}
\accepted{5 June 2004}

\usepackage{amsfonts, graphicx} 

\newtheorem{theorem}{Theorem}[section] 
\newtheorem{lemma}[theorem]{Lemma}     
\newtheorem{proposition}[theorem]{Proposition} 
\newtheorem{corollary}[theorem]{Corollary} 

\theoremstyle{definition}
\newtheorem{definition}[theorem]{Definition} 
\newtheorem{example}[theorem]{Example}       
  
\theoremstyle{remark}  
\newtheorem{remark}[theorem]{Remark}  
\newtheorem*{ack}{Acknowledgment}

\begin{document}
\title{Links associated with generic immersions of graphs}                    
\authors{Tomomi Kawamura}                  
\address{Department of Physics and Mathematics, 
Aoyama Gakuin University\\5-10-1, Fuchinobe Sagamihara, 
Kanagawa 229-8558, Japan }                  
\email{tomomi@gem.aoyama.ac.jp}

\begin{abstract}
As an extension of the class of algebraic links, A'Campo, Gibson, and
Ishikawa constructed links associated to immersed arcs and trees in a
two-dimensional disk.  By extending their arguments, we construct
links associated to immersed graphs in a disk, and show that such
links are quasipositive.
\end{abstract}

\primaryclass{57M25}                
\secondaryclass{57M27}              
\keywords{Divide, graph divide, 
quasipositive link, slice Euler characteristic, 
four-dimensional clasp number}                    

\maketitle 

\section{Introduction}\label{sect;intro}

In \cite{acampo1}, A'Campo constructed links of divides as an
extension of the class of algebraic links.  A {\em divide} is a
generic relative immersion of a disjoint union of arcs (and loops) in
a 2--dimensional disk.  In \cite{gibishi}, Gibson and Ishikawa
constructed links of free divides, non-relative immersions of
intervals in a 2--dimensional disk.  We review links of divides in
Section \ref{sect;revi}, and remark on links of free divides in
Section \ref{sect;defgd}.

In \cite{acampo}, A'Campo showed that any divide link is ambient
isotopic to a transverse $\mathbb{C}$--link.  A {\em transverse
$\mathbb{C}$--link} is the link represented as the transversal
intersection of an algebraic curve and the unit sphere in the
2--dimensional complex space $\mathbb{C}^2$ \cite{rutop90}.  An {\em
algebraic link} is the link of a singularity of an algebraic curve. An
algebraic link is a transverse $\mathbb{C}$--link; there exist
transverse $\mathbb{C}$--links which are not algebraic \cite{rualg}.

In \cite{ka6}, the author showed that links of divides and free
divides are quasipositive, by using the visualization algorithm due to
Hirasawa \cite{hira}.  A \emph{quasipositive braid} is a product of
braids which are conjugates of positive braids, and a
\emph{quasipositive link} is an oriented link which has a closed
quasipositive braid diagram.  A \emph{positive braid} is a product of
canonical generators of the braid group, that is, a braid which has a
diagram without negative crossings.  It is well known that any
algebraic link admits a representation as a closed positive braid.  In
\cite{rualg}, Rudolph showed that quasipositive links are transverse
$\mathbb{C}$--links.  In \cite{BO}, Boileau and Orevkov proved that
transverse $\mathbb{C}$--links are quasipositive, by using the theory
of pseudoholomorphic curves.

There exist quasipositive links which are not the links of any divides
or free divides \cite{ka6}.  For classification of such links, we are
interested in extending the class of links of divides and free
divides.  In \cite{gib;tree}, Gibson associated links with generic
immersions of trees in a 2--dimensional disk.  A {\em tree divide} is
such an immersion of trees.  He suggested to the author the
quasipositivity problem for such links.  In this paper, we construct
graph divide links as an extension of the class of links of divides or
free divides in Section \ref{sect;defgd}, show that tree divide links
constructed by Gibson can be represented as graph divide links in
Section \ref{sect;tree}, and prove that such links are quasipositive
in Section \ref{sect;qpos}.  Furthermore we determine some geometric
invariants for graph divide links, and show that there exists a
quasipositive link which is not a graph divide link.

\begin{ack}
The author was partially supported by JSPS Research Fellowships for
Young Scientists, while she was in Graduate School of Mathematical
Sciences, University of Tokyo.  She was also supported by Grant-in-Aid
for Young Scientists (B) (No.\ 15740044), MEXT.  

The author would like to thank Professor Masaharu Ishikawa and Doctor
William Gibson for their helpful suggestions and useful comments about
divides, free divides, tree divides, and graph divides.  She also
would like to thank them and Professor Toshitake Kohno for their
encouragement.  She would like to thank the referee for his/her
pointing out the ambiguity in the first manuscript.
\end{ack}

\section{Divide links and oriented divide links}\label{sect;revi}

In this section, we review links of divides defined by A'Campo
\cite{acampo1} and links of oriented divides defined by Gibson and
Ishikawa \cite{gibishiOD}.

By the argument due to A'Campo \cite{acampo1}, a link is associated to
any immersed arcs (and loops) in a disk as follows.

Let $D$ be a unit disk in the real plane $\mathbb{R}^2$, that is $D=\{
x=(x_1,x_2)\in \mathbb{R}^2 \mid |x|^2= x_1^2+x_2^2\leq 1\}$.  A
\emph{divide} $P$ is a generic relative immersion in the unit disk
$(D, \partial D)$ of a finite number of 1--manifolds, i.e., copies of
the unit interval $(I, \partial I)$ and the unit circle
\cite{acampo1,acampo,hira,ishi;Jpoly}.  We also call the image of such
an immersion a \emph{divide}.  A \emph{branch} of $P$ is each image of
the copies.  We shall call each image of the copies of the interval an
\emph{interval branch}, and each image of the copies of the circle a
\emph{circle branch}.

Let $T_xX$ be the tangent space 
at a point $x$ of a manifold $X$,  
and $TX$ be the tangent bundle over a manifold $X$. 
We identify the 3--sphere $S^3$ with the set 
\[ ST\mathbb{R}^2=\{ (x,u)\in T\mathbb{R}^2 \mid 
   x\in \mathbb{R}^2, u\in T_x\mathbb{R}^2, |x|^2+|u|^2=1 \} . \]
The \emph{link of a divide} $P$ is the set given by  
\[ L(P)=\{ (x,u)\in ST\mathbb{R}^2 \mid x\in P, u\in T_xP \} . \]
We orient the 3--sphere and the link $L(P)$ as follows.  
We identify the tangent bundle 
$T\mathbb{R}^2=\mathbb{R}^4$ with 
the 2--dimensional complex space $\mathbb{C}^2$ 
by the map 
\[ ((x_1,x_2),(u_1,u_2))
   \mapsto (x_1+\sqrt{-1}u_1,x_2+\sqrt{-1}u_2). \]
The tangent bundle $T\mathbb{R}^2$ is 
oriented by the complex orientation of $\mathbb{C}^2$, 
and the 3--sphere is naturally oriented 
by the complex orientation of the 4--ball 
\[ \{ (x,u)\in T\mathbb{R}^2 \mid 
   x\in \mathbb{R}^2, u\in T_x\mathbb{R}^2, 
   |x|^2+|u|^2\leq 1 \} . \]
Let $[a,b]$ be a small interval with $a<b$. 
Let $\phi \co  [a,b] \rightarrow D$ 
be an embedding whose image lies on $P$. 
We orient a part of the link $L(P)$ as the image of the map 
$t\mapsto (\phi (t), 
\displaystyle{\frac{\sqrt{1-|\phi (t)|^2}}{|\dot{\phi }(t)|}}
\dot{\phi }(t))$, 
where $\dot{\phi }(t)$ is the differential of $\phi (t)$. 
We can extend this orientation to $L(P)$. 
A {\em divide link} is an oriented link 
ambient isotopic to the link of some divide. 

\begin{figure}[ht!]
\begin{center}
\includegraphics[height=2.7cm]{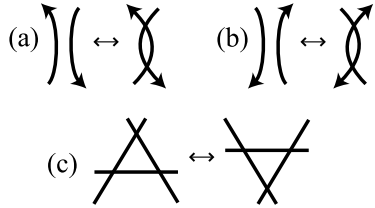}
\end{center}
\caption{Inverse self-tangency moves and a triangle move}
\label{fig;move}
\end{figure}

In \cite{gibishiOD}, 
Gibson and Ishikawa constructed 
links associated with oriented divides. 
An {\em oriented divide} is the image of a generic immersion of 
finite number of copies of the unit circle in the unit disk, 
with a specific orientation assigned to each immersed circle. 
The {\em link $L_{ori}(Q)$ of an oriented divide $Q$} 
is defined by 
\[ 
L_{ori}(Q)=\{ (x,u)\in ST\mathbb{R}^2 \mid 
x\in Q, u\in \vec{T}_xQ \}, 
\]
where $\vec{T}_xQ$ is the set of tangent vectors 
in the same direction as the assigned orientation of $Q$. 
The link $L_{ori}(Q)$ naturally inherits its orientation 
from $Q$. 
The ambient isotopy type of the link of an oriented divide 
does not change under 
the {\em inverse self-tangency moves} 
illustrated 
in Figure \ref{fig;move} (a) and (b), 
and the {\em triangle moves} 
illustrated 
in Figure \ref{fig;move} (c)
\cite{gibishi}. 
An {\em oriented divide link} is an oriented link 
ambient isotopic to the link of some oriented divide. 

\begin{figure}[ht!]
\begin{center}
\includegraphics[height=4cm]{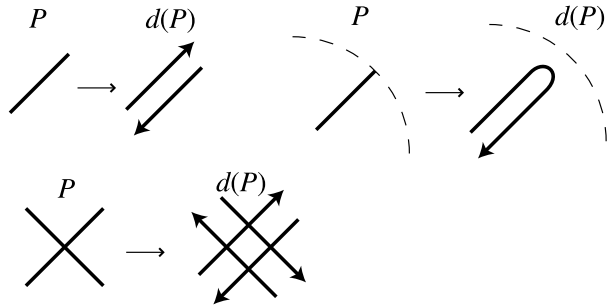}
\end{center}
\caption{Doubling method}
\label{fig;dp1}
\end{figure}

The link of a divide is isotopic to 
the link of the oriented divide 
obtained from the divide by the doubling method, 
which is the first step of 
the visualization algorithm due to Hirasawa \cite{hira}. 
Let $P$ be a divide. 
For each branch $B$ of $P$, 
we draw the boundary of a `very small' neighborhood of $B$ 
in the disk $D$, 
assigned with the clockwise orientation, 
as illustrated in Figure \ref{fig;dp1}, 
where interrupted curves represent $\partial D$. 
In particular, 
we draw a `sharp' around each double point of $P$, 
and draw a `hairpin curve' around each point of $\partial P$. 
We suppose that such hairpin curves lie in the interior of $D$. 
We denote by $d(P)$ the oriented divide 
obtained by the above algorithm. 

\begin{figure}[ht!]
\begin{center}
\includegraphics[height=10cm]{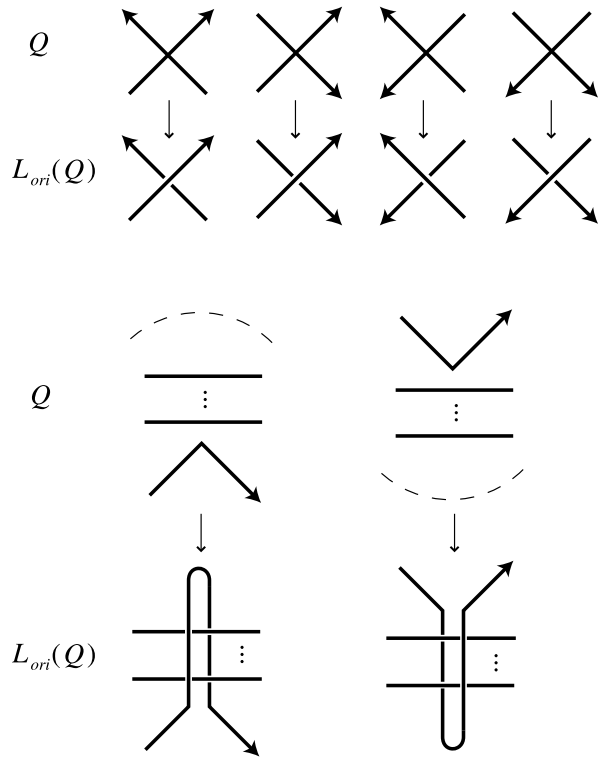}
\end{center}
\caption{Local link diagrams obtained from oriented divides}
\label{fig;diag1}
\end{figure}

We review the second step of 
the visualization algorithm due to Hirasawa \cite{hira}. 
We can apply his algorithm not only to divides 
but also to oriented divides. 
Let $Q$ be an oriented divide. 
A regular isotopy of $Q$ in the space of generic immersions 
does not change the ambient isotopy type of the link $L_{ori}(Q)$. 
Hence we may assume that 
oriented divides are linear with slope $\pm 1$ 
except near the `corners', 
where a branch quickly changes its slope from $\pm 1$ to $\mp 1$. 
We draw a link diagram for $Q$ 
as below: 

\begin{enumerate}
\item 
We replace each double point of $Q$ 
with the crossing 
as illustrated at the top of Figure \ref{fig;diag1}. 
\item 
For each $x_2$--maximal (resp.\ $x_2$--minimal) point  
whose tangent vector has same orientation as the $x_1$--axis, 
we change the diagram 
as illustrated at the bottom of Figure \ref{fig;diag1}, 
where interrupted curves represent $\partial D$, 
and horizontal arcs represent all arcs over (resp.\ under) 
the $x_2$--maximal (resp.\ $x_2$--minimal) point. 
\end{enumerate}

In \cite{hira}, 
Hirasawa showed that 
for any divide $P$, 
the diagram obtained by the above algorithm 
represents $L(P)$, the link of divide $P$. 
His argument implies that 
for any oriented divide $Q$, 
the diagram obtained by the second step of the above algorithm 
represents $L_{ori}(Q)$, the link of oriented divide $Q$. 
Therefore we obtain $L(P)=L_{ori}(d(P))$. 

\begin{remark}
In \cite{hira}, 
Hirasawa oriented the tangent bundle $TD$ by 
the coordinate $(x,u)=(x_1,x_2,u_1,u_2)$. 
Therefore the diagrams of links of divides in this paper are 
the mirror images of those in his paper, 
since the orientations of the 3--sphere are opposite. 
\end{remark}

\section{Links of graph divides}\label{sect;defgd}
%Def of links of graph divides. 
%In fact, we may consider only uni-trivalent graphs.  

In this section, 
we construct graph divide links, 
that is, 
links associated with 
generic immersions of finite graphs 
in a 2--dimensional disk $D$, 
as an extension of the class of 
divide links. 

%% Def of finite graph, degree of vertex  %%

A {\em graph divide} $P=(G, \varphi )$ is 
a generic immersion $\varphi \co  G\rightarrow D$ as follows 
or its image, 
where $G$ is a disjoint union of 
finite graphs and copies of the unit circle. 
Each graph may have loops and multiple edges. 
The singularities are only transversal double points 
of two arcs in edges and circles. 
We suppose that 
each point of $P\cap \partial D$
is the image of a vertex of degree $1$. 
We regard the unit interval as a finite graph. 
A \emph{branch} of $P$ is the image of each component of $G$. 
We shall call the image of an interval component
 an \emph{interval branch}, 
the image of a circle component 
a \emph{circle branch}, 
and the image of a tree component 
a \emph{tree branch}. 
The image of vertices of degree $1$ 
might not lie in the boundary of the unit disk. 
We call such an image a \emph{free endpoint} of $P$ 
and denote by $E_P$ the set of all free endpoints of $P$. 
We denote by $T_P$ the set of all vertices except 
free endpoints and  points in $\partial D$.
We denote $V_P=E_P\cup T_P$.  
If 
$T_P$ is empty and 
$\varphi$ is a non-relative immersion, 
$P$ is a {\em free divide} \cite{gibishi}. 
In \cite{gibishi}, Gibson and Ishikawa considered  
free divides with only interval branches, 
but we consider both of interval and circle branches 
in this paper. 
We note that a divide is also a free divide \cite{gibishi}, 
and that it is also a graph divide. 

We extend the definition of links of divides as follows. 
In the case of free divides, the argument is almost same as 
the visualized definition for links of free divide
due to Gibson and Ishikawa \cite{gibishi}. 
We give `signs' to vertices of 
a graph divide, 
because the link is not associated to a graph alone. 
Let $x$ be a vertex of $G$. 
We also denote the image $\varphi(x)$ by $x$. 
If $x$ lies in $\partial D$, 
$x$ does not need a sign. 
If $x$ is a point of $V_P$, 
we give $x$ a sign $\epsilon _x=+$ or $\epsilon _x=-$. 

\begin{figure}[ht!]
\begin{center}
\includegraphics[height=7cm]{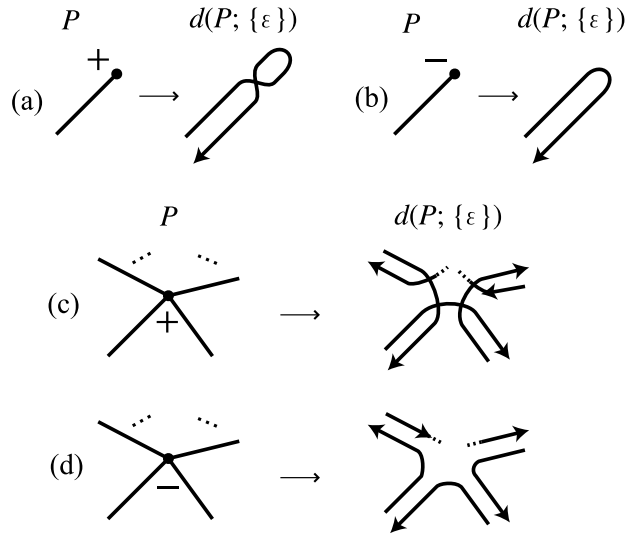}
\end{center}
\caption{The doubling of graph divide around vertices}
\label{fig;dp2}
\end{figure}

For a given graph divide $P=(G, \varphi)$ 
and given signs of vertices, 
we construct an oriented divide 
$d(P;\{ \epsilon _x\} _{x\in V_P})$
by extending a doubling method as follows. 
For each branch $B$ of $P$, 
we draw the boundary of `very small' neighborhood of $B$ 
in the disk $D$, assigned with the clockwise orientation, 
in the same way as that for divides, 
except near $x\in V_P$. 
Around $x\in E_P$ with $\epsilon _x=-$, 
we draw a `hairpin curve', 
as illustrated in 
Figure \ref{fig;dp2} (b). 
Around $x\in E_P$ with $\epsilon _x=+$, 
we draw a `loop', 
as illustrated in 
Figure \ref{fig;dp2} (a). 
Around $x\in T_P$ with $\epsilon _x=-$, 
we draw oriented curves 
such that 
each curve approaches $x$ along an edge and 
turns to its neighbor edge on the left, 
as illustrated in 
Figure \ref{fig;dp2} (d). 
Around $x\in T_P$ with $\epsilon _x=+$, 
we draw oriented curves 
such that 
each curve approaches $x$ along an edge and 
turns to its neighbor edge on the right, 
as illustrated in 
Figure \ref{fig;dp2} (c). 
We denoted the obtained curves by 
$d(P;\{ \epsilon _x\} _{x\in V_P})$
and call it 
the {\em doubling} of graph divide $P$ 
with signs 
$\{ \epsilon _x\} _{x\in V_P}$. 

\begin{definition}
The {\em link of a graph divide} $P$ 
is the set given by
\[ 
L(P;\{ \epsilon _x\} _{x\in V_P})
=L_{ori}(
d(P;\{ \epsilon _x\} _{x\in V_P})
). \]
\end{definition}

We note that 
the link of a given graph divide 
depends on signs of vertices. 
For fixed signs of vertices  
$(\{ \epsilon _x\} _{x\in V_P})$, 
a regular isotopy of $P$ in the space of generic immersions 
does not change the isotopy type of the oriented divide 
$d(P;\{ \epsilon _x\} _{x\in V_P})$, 
therefore it does not change the ambient isotopy type of 
the link 
$L(P;\{ \epsilon _x\} _{x\in V_P})$. 
A \emph{graph divide link} is the oriented link 
ambient isotopic to 
the link of some graph divide 
with some signs of vertices. 

Furthermore some transformations on a graph divide 
do not change the isotopy type of 
the link. 

\begin{figure}[ht!]
\begin{center}
\includegraphics[height=4cm]{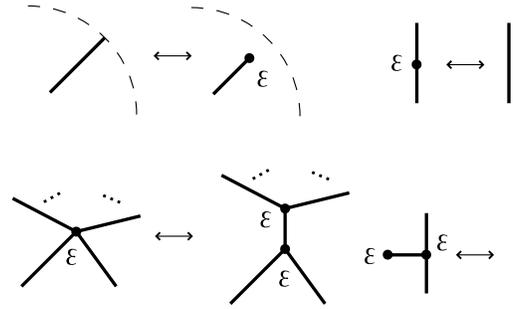}
\end{center}
\caption{Transformations on a graph divide 
not changing the link type}
\label{fig;trans}
\end{figure}

\begin{lemma}\label{lem;trans}
The transformations on a graph divide 
illustrated in Figure \ref{fig;trans} 
do not change the ambient isotopy type of the link. 
\end{lemma}

\begin{proof}
For each transformation illustrated in Figure \ref{fig;trans}, 
we consider the doubling of the graph divides 
illustrated in Figure \ref{fig;dp2}.  
They are changed to each other 
by diffeomorphisms of $D$, 
triangle moves and inverse self-tangency moves. 
\end{proof}

By means of Lemma \ref{lem;trans}, 
for any graph divide $P=(G,\varphi )$, 
there exists a graph divide $P'=(G',\varphi ')$,
where $G'$ is an union of 
uni-trivalent graphs and copies of circles,   
such that $L(P';\{ \epsilon _x\} _{x\in V_{P'}})$ is 
ambient isotopic 
to $L(P;\{ \epsilon _x\} _{x\in V_P})$. 

%\begin{remark}
The above definition of a graph divide 
is a natural extension of a (free) divide link. 
If $E_P=V_P$ holds, the link $L(P;\{ \epsilon _x\} _{x\in V_P})$ 
is ambient isotopic to the link of a free divide 
defined by Gibson and Ishikawa \cite{gibishi}. 
%\end{remark}

\begin{figure}[ht!]
\begin{center}
\includegraphics[height=7.5cm]{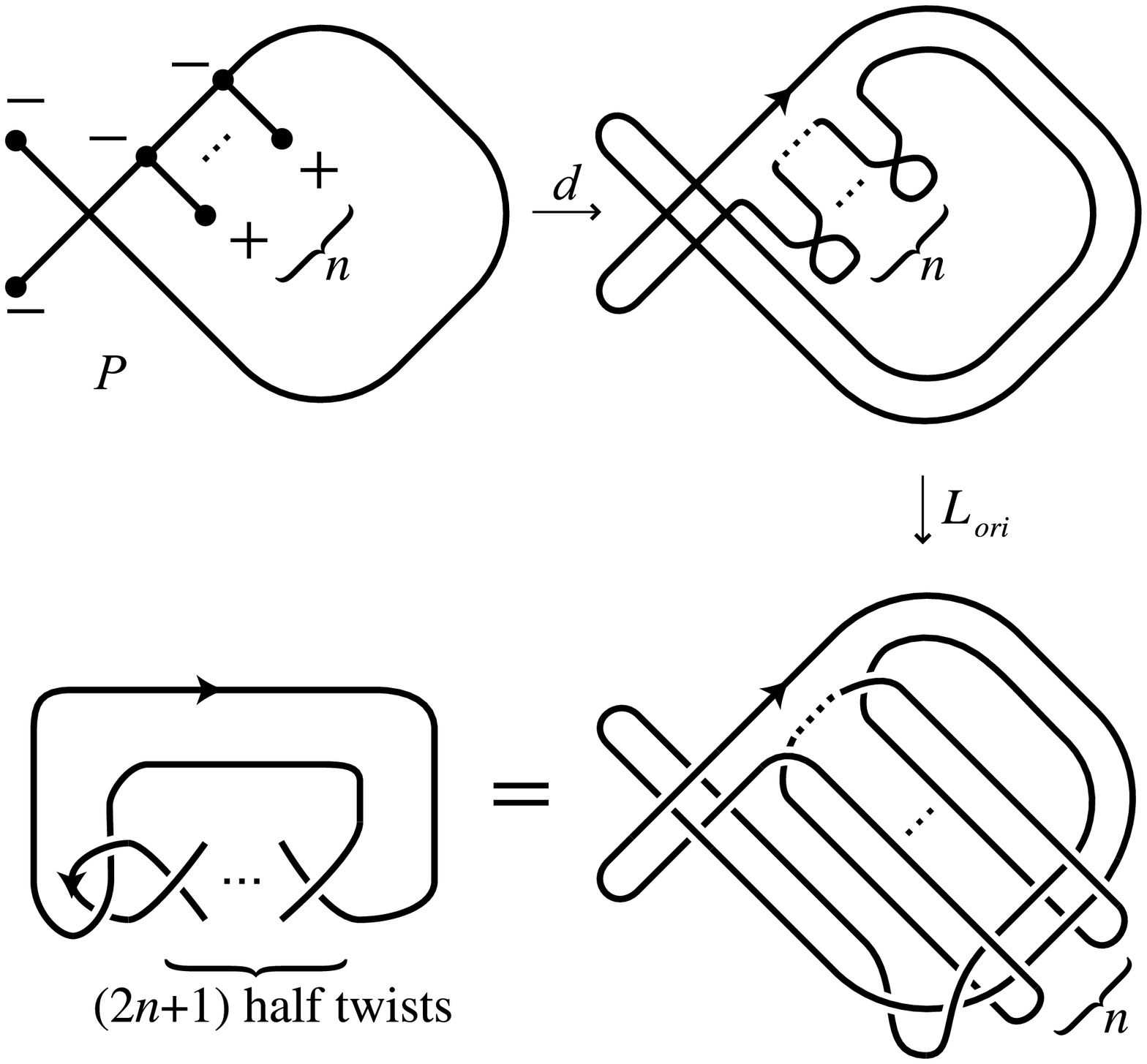}
\end{center}
\caption{An example of a (non-fibered) graph divide knot}
\label{fig;fish1}
\end{figure}

\begin{example}(Cf.\ Gibson \cite{gib;tree})\label{ex;fish1}\qua
For a graph divide $P$ with signed vertices illustrated 
in Figure \ref{fig;fish1}, 
the doubling of $P$ 
is illustrated as 
the right of $P$.  
Then 
the link of $P$ 
is the knot 
illustrated at the bottom of Figure \ref{fig;fish1}. 
It is known that this knot is not fibered 
if $n$ is a positive integer. 
Then it is not a divide link since divide links are all fibered. 
\end{example}

\begin{figure}[ht!]
\begin{center}
\includegraphics[height=6.7cm]{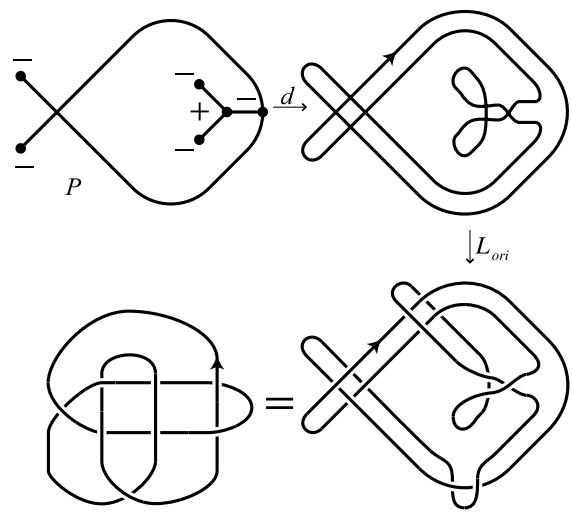}
\end{center}
\caption{An example of a graph divide knot}
\label{fig;fish2}
\end{figure}

\begin{example}\label{ex;fish2}
For a graph divide $P$ with signed vertices 
illustrated 
in Figure \ref{fig;fish2}, 
the doubling of $P$ 
is illustrated as 
the right of $P$. 
Then 
the link of $P$ is the knot 
illustrated at the bottom of Figure \ref{fig;fish2}. 
This knot is the mirror image of $8_{21}$ 
in the table of Rolfsen \cite{rolfsen}. 
It is well known that the knot $8_{21}$ is fibered. 
\end{example}

In the next section, %Section \ref{sect;tree}, 
we compare the above definition of the link of a graph divide 
with the construction of the link of a tree divide 
defined by Gibson in \cite{gib;tree}.  
We use the following lemma, 
where 
we denote by $\overline{Q}$ 
the oriented divide obtained from 
a given oriented divide $Q$ 
by reversing the orientations of all branches. 
Gibson and Ishikawa showed as Proposition 3.1 in \cite{gibishiOD} 
that 
$L_{ori}(\overline{Q})$ is the same link as $L_{ori}(Q)$ 
but with the opposite orientations on all components of $L_{ori}(Q)$. 

\begin{lemma}\label{lem;oppositeori}
Let $P=(G_1\sqcup G_2, \varphi )$ be a graph divide and 
we denote $(G_j, \varphi |_{G_j})$ by $P_j=(G_j, \varphi _j)$ 
for $j=1,2$. 
We give each vertex $x\in V_P=V_{P_1}\cup V_{P_2}$ 
a sign $\epsilon _x$. 
Then the link 
$\displaystyle L(P;
\{ \epsilon _x\} _{x\in V_{P_1}}\cup \{ -\epsilon _x\} _{x\in V_{P_2}})$ 
is ambient isotopic to 
the link $L_{ori}(
d(P_1;\{ \epsilon _x\} _{x\in V_{P_1}})\cup 
\overline{d(P_2;\{ \epsilon _x\} _{x\in V_{P_2}})})$, 
where $\{ -\epsilon _x\} _{x\in V_{P_2}}$ 
is the set of the signs of $x\in V_{P_2}$ defined by $-\epsilon _x$. 
\end{lemma}

\begin{figure}[ht!]
\begin{center}
\includegraphics[height=9cm]{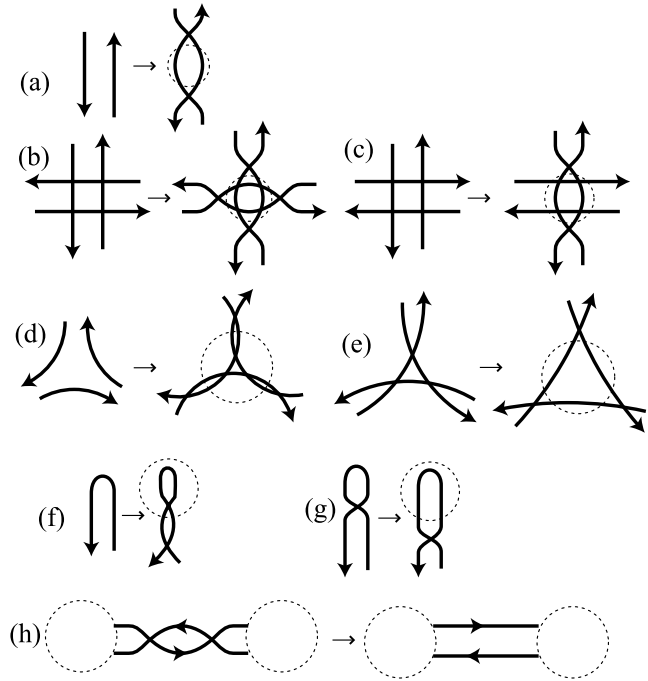}
\end{center}
\caption{Transformations on 
$\overline{d(P_2;\{ \epsilon _x\} _{x\in V_{P_2}})}$}
\label{fig;opori}
\end{figure}

\begin{proof}
By Lemma \ref{lem;trans}, 
we may suppose that $P\cap \partial D=\emptyset$ 
and the degree of each vertex of $G$ is $1$ or $3$. 
We transform each part of 
$\overline{d(P_2;\{ \epsilon _x\} _{x\in V_{P_2}})}$ 
as illustrated in Figure \ref{fig;opori} 
(a), (b), (c), (d), and (f)
by inverse self-tangency moves and triangle moves,  
and as illustrated in Figure \ref{fig;opori} (e) and (g) 
by diffeomorphisms of $D$. 
We denote by $Q_2$ this new oriented divide obtained 
from $\overline{d(P_2;\{ \epsilon _x\} _{x\in V_{P_2}})}$. 
Each of the parts in 
the small disks with dotted boundary  
in Figure \ref{fig;opori} 
is same as the assigned part of 
the doubling of $d(P_2;\{ -\epsilon _x\} _{x\in V_{P_2}})$.  
We remove the double points of $Q_2$ 
in the exterior of such disks 
by inverse self-tangency moves 
as illustrated in Figure \ref{fig;opori} (h). 
The finally obtained oriented divide is 
same as $d(P_2;\{ -\epsilon _x\} _{x\in V_{P_2}})$. 
Then the oriented divide 
$d(P_1;\{ \epsilon _x\} _{x\in V_{P_1}})\cup 
\overline{d(P_2;\{ \epsilon _x\} _{x\in V_{P_2}})}$ 
is changed to 
$d(P;
\{ \epsilon _x\} _{x\in V_{P_1}}\cup \{ -\epsilon _x\} _{x\in V_{P_2}})$ 
by diffeomorphisms of $D$, 
inverse self-tangency moves, and triangle moves. 
Therefore the link $L_{ori}(
d(P_1;\{ \epsilon _x\} _{x\in V_{P_1}})\cup 
\overline{d(P_2;\{ \epsilon _x\} _{x\in V_{P_2}})})$ 
is ambient isotopic to 
$\displaystyle L(P;
\{ \epsilon _x\} _{x\in V_{P_1}}\cup \{ -\epsilon _x\} _{x\in V_{P_2}})$. 
\end{proof}

\begin{remark}
Let $P=(I,\varphi )$ be a free divide with a single interval branch 
and $V_P=\{ x_1, x_2\}$. 
In \cite{gibishi} Gibson and Ishikawa said that 
the \emph{sign} of $P$ is \emph{even} 
if $\epsilon _{x_1}=\epsilon _{x_2}$ holds, 
and otherwise \emph{odd}. 
They showed that 
the isotopy class of the knot of such a free divide $P$ 
depends only on the sign of $P$. 
Actually, 
the orientation of $L(P;\{ -\epsilon _{x_1}, -\epsilon _{x_2})\} )$ 
is the reverse of $L(P;\{ \epsilon _{x_1}, \epsilon _{x_2})\} )$. 
These facts are included in Lemma \ref{lem;oppositeori}. 
\end{remark}

\section{A relation with the construction due to Gibson
%tree divides
}\label{sect;tree}
%Links of tree divides defined by Gibson 
%are also links of graph divides.  

We call a graph divide $P=(G, \varphi )$ 
a {\em tree divide} 
if $G$ is an union of trees \cite{gib;tree}. 
In \cite{gib;tree} Gibson defined 
the link of a tree divide 
by the different argument from 
that of the link of a graph divide in this paper. 
In this section, we review his definition 
and show that 
tree divide links defined by him 
are represented as 
graph divide links defined in 
Section \ref{sect;defgd}.

\begin{figure}[ht!]
\begin{center}
\includegraphics[height=11cm]{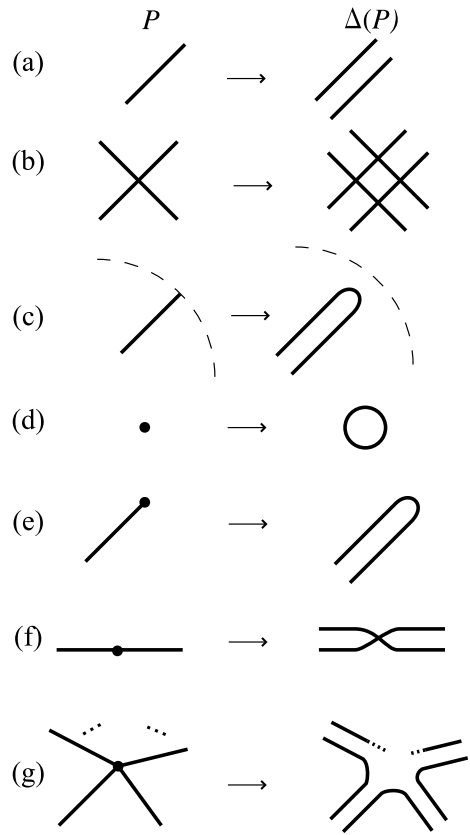}
\end{center}
\caption{The doubling $\Delta (P)$ of a tree divide $P$}
\label{fig;wtree}
\end{figure}

For a given tree divide $P=(G, \varphi)$, 
Gibson defined a {\em doubling} of $P$ 
by the following argument. 
We assume that an union of trees 
$G$ may have isolated vertices. 
%We may suppose that 
%the tree divide $P=(G, \varphi )$ does not intersect
%the boundary of the 2--disk $D$. 
We construct a new non-oriented divide $\Delta (P)$ as follows. 
For each branch $B$ of $P$, 
we draw the boundary of `very small' neighborhood of $B$ 
in the disk $D$ except near $x\in V_P$ 
as 
illustrated 
in Figure \ref{fig;wtree} (a), (b), and (c), 
in the almost same way as that for the doubling in Section \ref{sect;revi}. 
In this step we do not give any orientation. 
Around each vertex we connect these curves 
as illustrated in Figure \ref{fig;wtree} as follows. 
Around each isolated vertex, 
we draw a small circle as illustrated in Figure \ref{fig;wtree} (d). 
Around $x\in E_P$, 
we draw a `hairpin curve' as illustrated in Figure \ref{fig;wtree} (e). 
Around each vertex of degree $2$, 
we draw two curves which cross each other once transversely 
as illustrated in Figure \ref{fig;wtree} (f). 
Around each of the remaining vertices, 
as illustrated in Figure \ref{fig;wtree} (g), 
we draw non-oriented curves 
the same way as the doubling around the vertex signed `$-$' 
in the definition of a doubling of graph divide 
in Section \ref{sect;defgd}. 

We give $\Delta (P)$ an arbitrary orientation $o$. 
Then we obtain an oriented divide $(\Delta (P), o)$. 
Gibson defined 
the {\em link of a tree divide} $P$
as the link $L_{ori}(\Delta (P), o)$. 

Applying the above argument 
to an immersion of the graphs 
which are not trees 
is not a natural extension of 
the original definition of the link of a divide, 
because if we do that, 
the doubling of some regular arc of $P$ 
might induce two parallel arcs with same orientation 
after we orient $\Delta (P)$. 
By the following result, 
we see that 
our definition of links of graph divides is 
a natural extension of links of divides and tree divides. 

\begin{proposition}
For any tree divide $P$, 
the link $L_{ori}(\Delta (P), o)$ can be 
represented as 
a graph divide link defined in Section \ref{sect;defgd}. 
\end{proposition} 

\begin{figure}[ht!]
\begin{center}
\includegraphics[height=6.5cm]{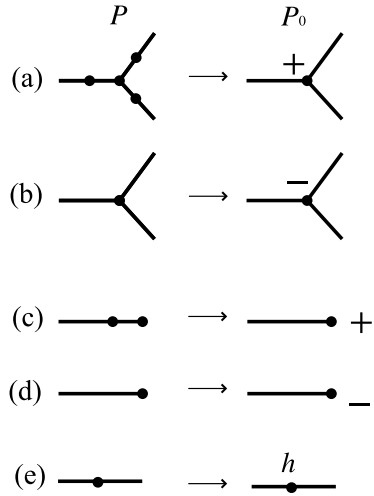}
\end{center}
\caption{Transformations from $P$ to $P_0$}
\label{fig;trgr}
\end{figure}

\begin{figure}[ht!]
\begin{center}
\includegraphics[height=4.5cm]{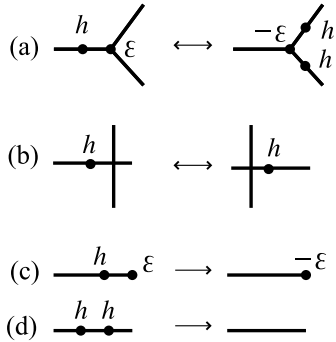}
\end{center}
\caption{The algorithm to remove vertices of degree 2}
\label{fig;trgr2}
\end{figure}

\begin{proof} Let $P=(G, \varphi )$ be a tree divide. 
By definition of the link of a tree divide, 
each isolated vertex may be replaced with a small embedded interval, 
and we may suppose that 
$P$ does not intersect $\partial D$ and  
the degree of each vertex of $G$ is $1$, $2$, or $3$. 
We construct a new tree divide $P_0$ and 
give 
the sign $\epsilon =+,-$ or the symbol $h$ 
to each vertex 
by replacing neighborhood of vertices 
as follows:

\begin{enumerate}
\item 
We give the sign `$+$' to each vertex of degree 3 
adjacent to three vertices of degree 2, 
and remove these vertices of degree 2, 
as illustrated in Figure \ref{fig;trgr} (a). 
We give the sign `$-$' 
to each of the other vertices of degree 3, 
as illustrated in Figure \ref{fig;trgr} (b). 
\item 
We give the sign $+$ to each endpoint 
adjacent to a vertex of degree 2, 
and remove this vertex of degree 2, 
as illustrated in Figure \ref{fig;trgr} (c). 
We give the sign `$-$' 
to each of the other endpoints, 
as illustrated in Figure \ref{fig;trgr} (d). 
\item 
After applying the above two steps, 
we give the symbol 
`$h$' 
to each of the other vertices of degree 2, 
in order to distinguish 
the doubling around a vertex of degree 2 
in this section 
from that defined in Section \ref{sect;defgd}.  
\end{enumerate}

We define the doubling of $P_0$, 
$\Delta (P_0; \{ \eta _x\} _{x\in V_{P_0}})$ as $\Delta (P)$, 
where $\eta _x$ is 
$h$ or $+$ or $-$ defined for each $x\in V_{P_0}$ as above. 
The moves (c) and (d) in Figure \ref{fig;trgr2} reduce 
the number of the vertices of degree $2$, 
changing the sign of the endpoint. 
The moves (a) and (b) in Figure \ref{fig;trgr2} can 
bring the vertices of degree $2$ close to the endpoints, 
changing the sign of the vertex of degree 3. 
By repeating these moves finitely many times, 
We obtain a tree divide $P_1=(G_1, \varphi)$ 
which has no vertex of degree $2$, 
and obtain signs of vertices $\{ \epsilon _x\} _{x\in V_{P_1}}$. 

Each of the transformations 
illustrated in Figure \ref{fig;trgr2} 
does not change 
the doubling of tree divides 
up to inverse self-tangency moves 
and triangle moves. 
Then, for any orientation $o$ of $\Delta (P)$, 
the oriented divide $(\Delta (P), o)$ can be transformed to 
$\Delta (P_1; \{ \epsilon _x\} _{x\in V_{P_1}})$ 
with some orientation $o_1$  
by inverse self-tangency moves and triangle moves. 
The oriented divide $d(P_1;\{ \epsilon _x\} _{x\in V_{P_1}})$ 
is also obtained as 
$\Delta (P_1; \{ \epsilon _x\} _{x\in V_{P_1}})$ 
with some orientation $o'_1$. 
For each branch $B$ of $P_1$ with 
$\displaystyle o_1|_{\Delta (B; \{ \epsilon _x\} _{x\in V_B})}
\neq o'_1|_{\Delta (B; \{ \epsilon _x\} _{x\in V_B})}$, 
we change the signs of all vertices of $B$. 
We denote new signs by $\{ \epsilon '_x\} _{x\in V_{P_1}}$. 
By Lemma \ref{lem;oppositeori}, 
the link $L(P_1;\{ \epsilon '_x\} _{x\in V_{P_1}})$ is 
ambient isotopic to 
$L_{ori}(\Delta (P_1;\{ \epsilon _x\} _{x\in V_{P_1}}), o_1)$. 
Therefore it is ambient isotopic to 
the link $L_{ori}(\Delta (P), o)$. 
\end{proof}

\section{Quasipositivity}\label{sect;qpos}
%Quasipositivity of graph divide links.

By using the almost same argument as 
that due to Ishikawa \cite{ishi;Jpoly},
we have the following theorem. 

\begin{theorem}\label{thm;grdivqp}
Links of graph divides are quasipositive. 
\end{theorem}

\begin{figure}[ht!]
\begin{center}
\includegraphics[height=5.5cm]{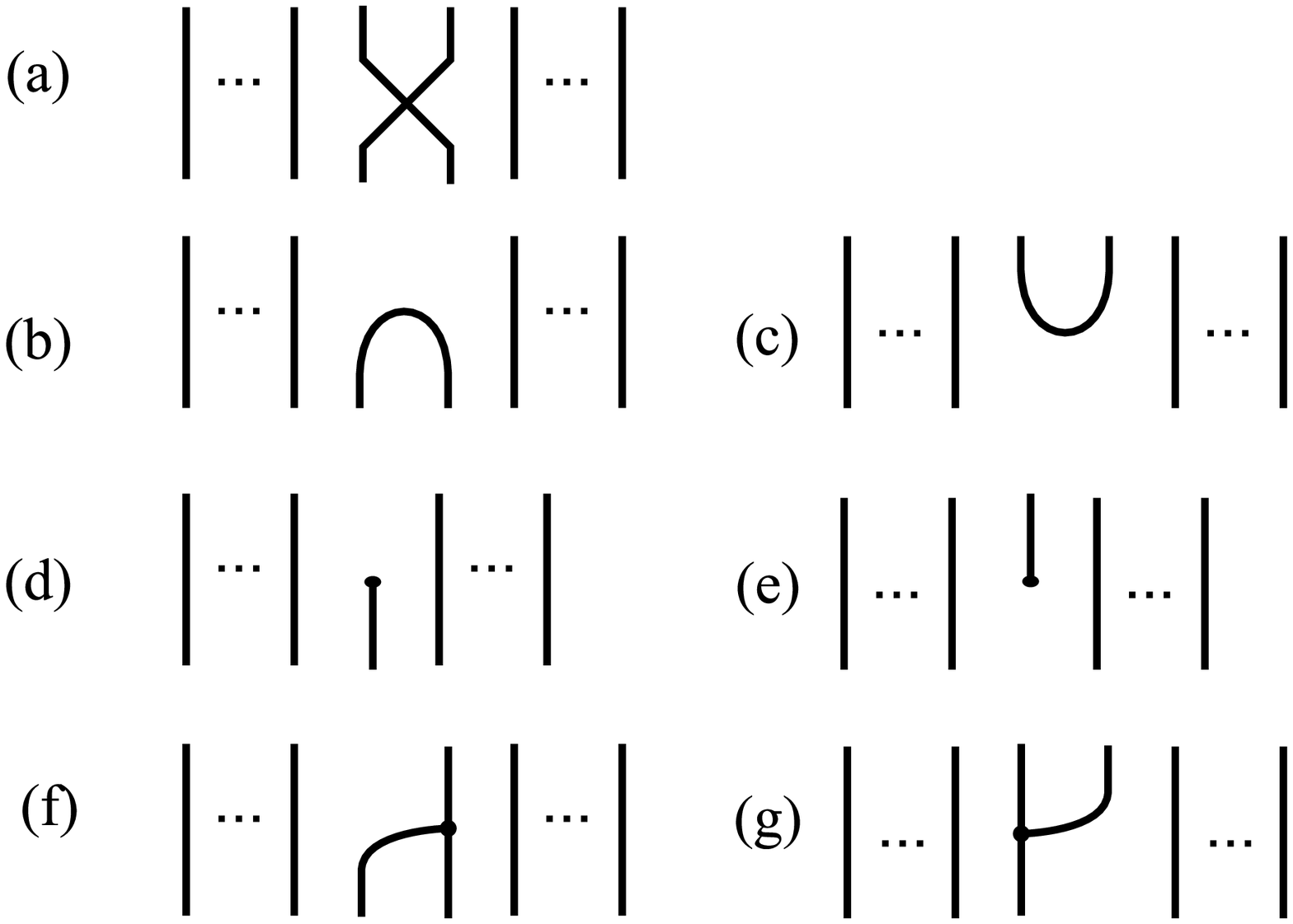}
\end{center}
\caption{Tangles}
\label{fig;tangle}
\end{figure}

In the proof of this theorem, 
we represent a graph divide as a tangle product 
constructed as follows.  
We restate a tangle product 
which Ishikawa defined for a divide in \cite{ishi;Jpoly}. 
In this paper, we suppose that 
a \emph{tangle} consists of some or no vertical lines 
and each of the following parts 
as illustrated in Figure \ref{fig;tangle}: 

\begin{enumerate}
\item a pair of crossed curves illustrated in Figure \ref{fig;tangle} (a),  
\item a folding curve including an $x_2$--maximal point 
  illustrated in Figure \ref{fig;tangle} (b),  
\item a folding curve including an $x_2$--minimal point 
  illustrated in Figure \ref{fig;tangle} (c), 
\item a vertical line with a vertex of degree $1$ 
  illustrated in Figure \ref{fig;tangle} (d) and (e), 
\item a vertical line and a curve with a vertex of degree $3$ 
  illustrated in Figure \ref{fig;tangle} (f) and (g).  
\end{enumerate}

A \emph{tangle product} is a product of these tangle with 
well-defined connections. 
The graph divide $P$ in Figure \ref{fig;exFP} 
is a tangle product representation 
for the graph divide of Example \ref{ex;fish2} 
illustrated in Figure \ref{fig;fish2}. 

\begin{remark}
In \cite{ishi;Jpoly} Ishikawa constructed a tangle 
as horizontal lines instead of vertical lines. 
In his construction, some folding curves have marks 
in order to establish the Kauffman state model on divides. 
Here we do not mention about a mark 
since we do not need it in the present article. 
\end{remark}

\begin{proof}[Proof of Theorem \ref{thm;grdivqp}]
Let $P=(G,\varphi )$ be a graph divide. 
By Lemma \ref{lem;trans}, 
we may suppose that for each vertex of $G$, 
the degree is equal to $1$ or $3$. 
A regular isotopy of $P$ in the space of generic immersions 
does not change the isotopy type of the link $L(P)$. 
Hence, 
we may assume that 
$P$ is represented as a tangle product. 
The horizontal arcs illustrated 
in Figure \ref{fig;FPe1}, \ref{fig;FPe2}, and \ref{fig;FPe3} 
represent all parts over or under 
each of vertices, and $x_2$--maximal or minimal points. 
Furthermore, we may assume that for any $a\in [-1,1]$
the number of arc components of $(p_1|_P)^{-1}(a)$ 
is less than $2$ except when 
they are connected by the neighborhood of double points in $P$, 
where $p_1$ is 
the projection map from the disk D to the $x_1$--axis. 

For a given graph divide $P$ and signs, 
we construct 
an immersed 2--manifold in the 3--sphere, 
$F(P;\{ \epsilon _x\} _{x\in V_P})$, 
which consists of disks and bands 
as determined by the following algorithm. 

\begin{figure}[ht!]
\begin{center}
\includegraphics[height=8.5cm]{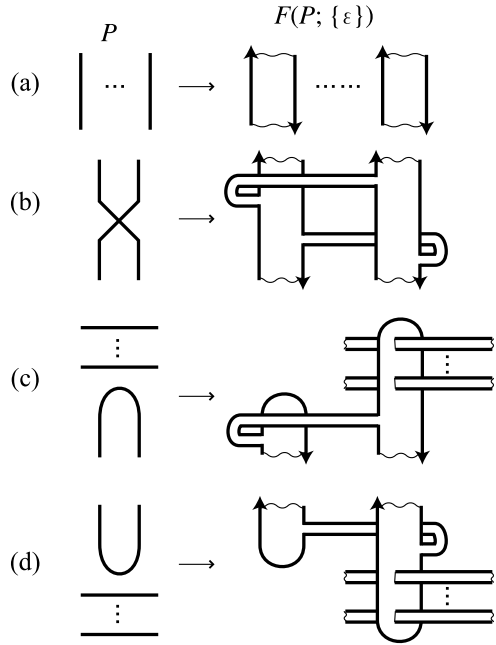}
\end{center}
\caption{Steps 1, 2, and 3 of the construction of 
$F(P;\{ \epsilon _x\} _{x\in V_P})$}
\label{fig;FPe1}
\end{figure}

\begin{enumerate}
\item 
For each vertical arc $\alpha$ in $P$, 
we construct a wide band along $\alpha$ in the disk $D$
as illustrated 
in Figure \ref{fig;FPe1} (a). 
This band is a part of `disks' of 
$F(P;\{ \epsilon _x\} _{x\in V_P})$, 
that is, 
each of such disks is 
the natural connected sum of the neighborhood of 
1--dimensional components of $(p_1|_P)^{-1}(a)$ in $D$ 
for some $a$. 
\item 
Around each double point of $P$, 
we construct two parts of disks and 
two narrow bands connecting them  
as illustrated in Figure \ref{fig;FPe1} (b). 
\item 
Around each of 
$x_2$--maximal points and $x_2$--minimal points of $P$, 
we construct two parts of disks and 
a narrow band connecting them, 
as illustrated in Figure \ref{fig;FPe1} (c) and (d). 
We construct intersections of a disk and 
the narrow bands corresponding to all small arcs in $P$ 
over the $x_2$--maximal point or 
under the $x_2$--minimal point.  
\item 
Around each free endpoint of $P$, 
we construct a part of a disk 
and intersections of it and
narrow bands corresponding to all small arcs in $P$ 
over or under the endpoint 
as illustrated in Figure \ref{fig;FPe2}. 
\item 
Around each vertex of degree $3$ of $P$, 
we construct two parts of disks and 
a narrow band connecting them, 
as illustrated in Figure \ref{fig;FPe3}. 
We construct intersections of a disk and 
the narrow bands corresponding to all small arcs in $P$ 
over or under the vertex. 
\end{enumerate}

We note that 
the narrow bands constructed in Step 2, 3, and 5 
may intersect some disks 
as explained in Step 3, 4, and 5. 

\begin{figure}[ht!]
\begin{center}
\includegraphics[height=3.75cm]{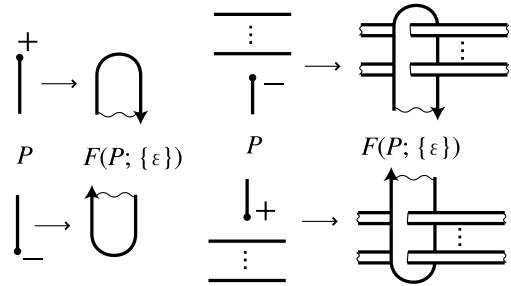}
\end{center}
\caption{Step 4 of the construction of $F(P;\{ \epsilon _x\} _{x\in V_P})$}
\label{fig;FPe2}
\end{figure}

\begin{figure}[ht!]
\begin{center}
\includegraphics[height=9.75cm]{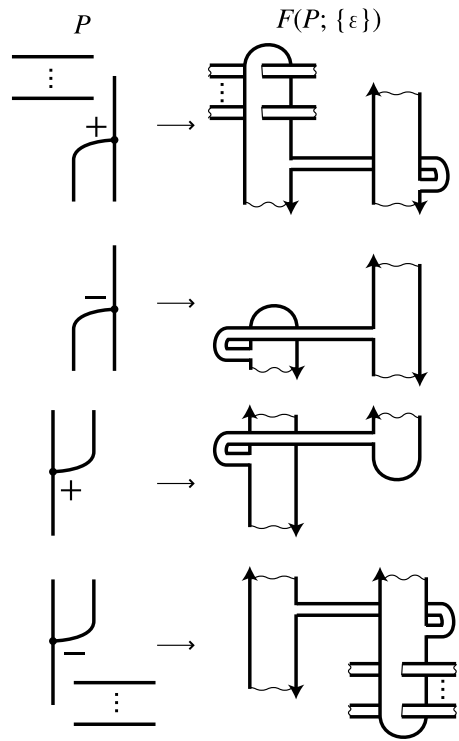}
\end{center}
\caption{Step 5 of construction of $F(P;\{ \epsilon _x\} _{x\in V_P})$}
\label{fig;FPe3}
\end{figure}

\begin{figure}[ht!]
\begin{center}
\includegraphics[height=5cm]{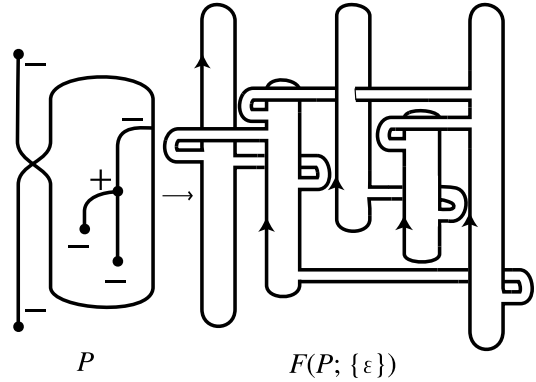}
\end{center}
\caption{An example of $F(P;\{ \epsilon _x\} _{x\in V_P})$}
\label{fig;exFP}
\end{figure}

Figure \ref{fig;exFP} is an example of 
the pair of $P$ and $F(P;\{ \epsilon _x\} _{x\in V_P})$. 
We suppose that 
the boundary of the above immersed 2--manifold, 
that is,  
$\partial F(P;\{ \epsilon _x\} _{x\in V_P})$ 
is oriented clockwise in the diagram obtained 
by the above algorithm 
as illustrated in Figure \ref{fig;FPe1}, \ref{fig;FPe2}, and \ref{fig;FPe3}. 
We regard the diagram of
$\partial F(P;\{ \epsilon _x\} _{x\in V_P})$ 
as a closed braid diagram. 
Each narrow band corresponds to 
a \emph{quasipositive band}, 
a conjugate braid with a canonical generator of the braid group. 
Then $\partial F(P;\{ \epsilon _x\} _{x\in V_P})$ 
is quasipositive.  
In the case of Figure \ref{fig;exFP}, 
the link $\partial F(P;\{ \epsilon _x\} _{x\in V_P})$ 
is the closure of the quasipositive braid 
\begin{eqnarray*}
&\sigma _1\sigma _4(\sigma _4^{-1}\sigma _3\sigma _2\sigma _3^{-1}\sigma _4)
\sigma _1\sigma _3(\sigma _4^{-1}\sigma _3^{-1}\sigma _2\sigma _3\sigma _4)\\
=&\sigma _1\sigma _4(\sigma _2^{-1}\sigma _3\sigma _4\sigma _3^{-1}\sigma _2)
\sigma _1\sigma _3(\sigma _2\sigma _3\sigma _4\sigma _3^{-1}\sigma _2^{-1}),
\end{eqnarray*} 
where $\sigma _1$, $\sigma _2$, $\sigma _3$, and $\sigma _4$ are 
the canonical generators of the 5--braid group. 
Comparing with the definition of the link of a graph divide, 
$\partial F(P;\{ \epsilon _x\} _{x\in V_P})$ 
is ambient isotopic to 
$L(P;\{ \epsilon _x\} _{x\in V_P})$. 
Therefore the link $L(P;\{ \epsilon _x\} _{x\in V_P})$ is quasipositive. 
\end{proof}

The {\em braid index} of a link $L$ is 
the minimal number of strings required to represent $L$ 
as a closed braid. 
The above argument implies the following result. 

\begin{proposition}
Let $P$ is a graph divide represented as 
a tangle product in the proof of Theorem \ref{thm;grdivqp}. 
Let $v$ be the number of vertices of $P$ 
and $m$ be the number of $x_2$--maximal points and 
$x_2$--minimal points of $P$. 
Then the braid index of the link of $P$ is 
not greater than 
$\displaystyle \frac{1}{2}(v+2m)$. 
\qed
\end{proposition}

\begin{figure}[ht!]
\begin{center}
\includegraphics[height=5cm]{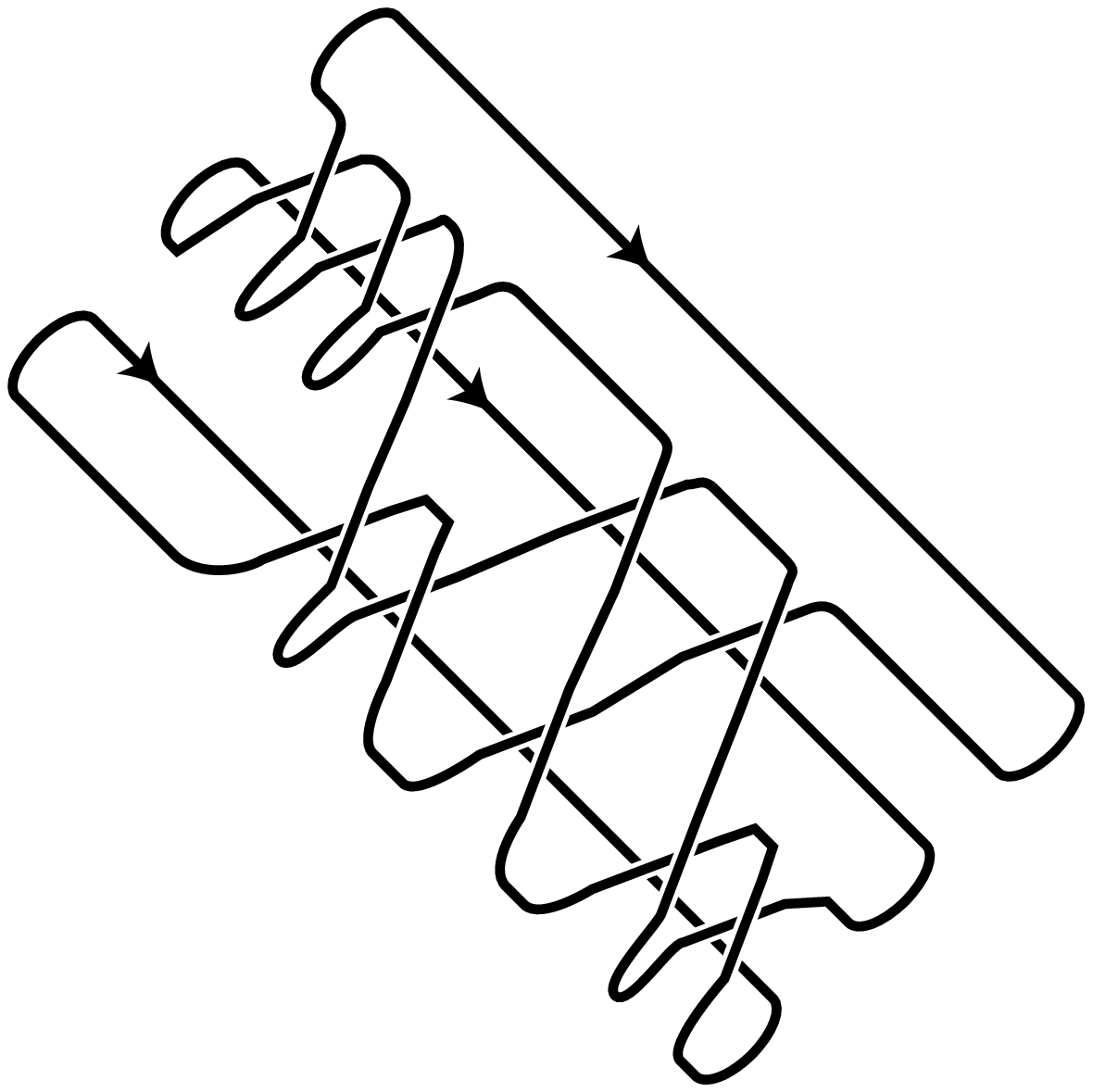}
\includegraphics[height=5cm]{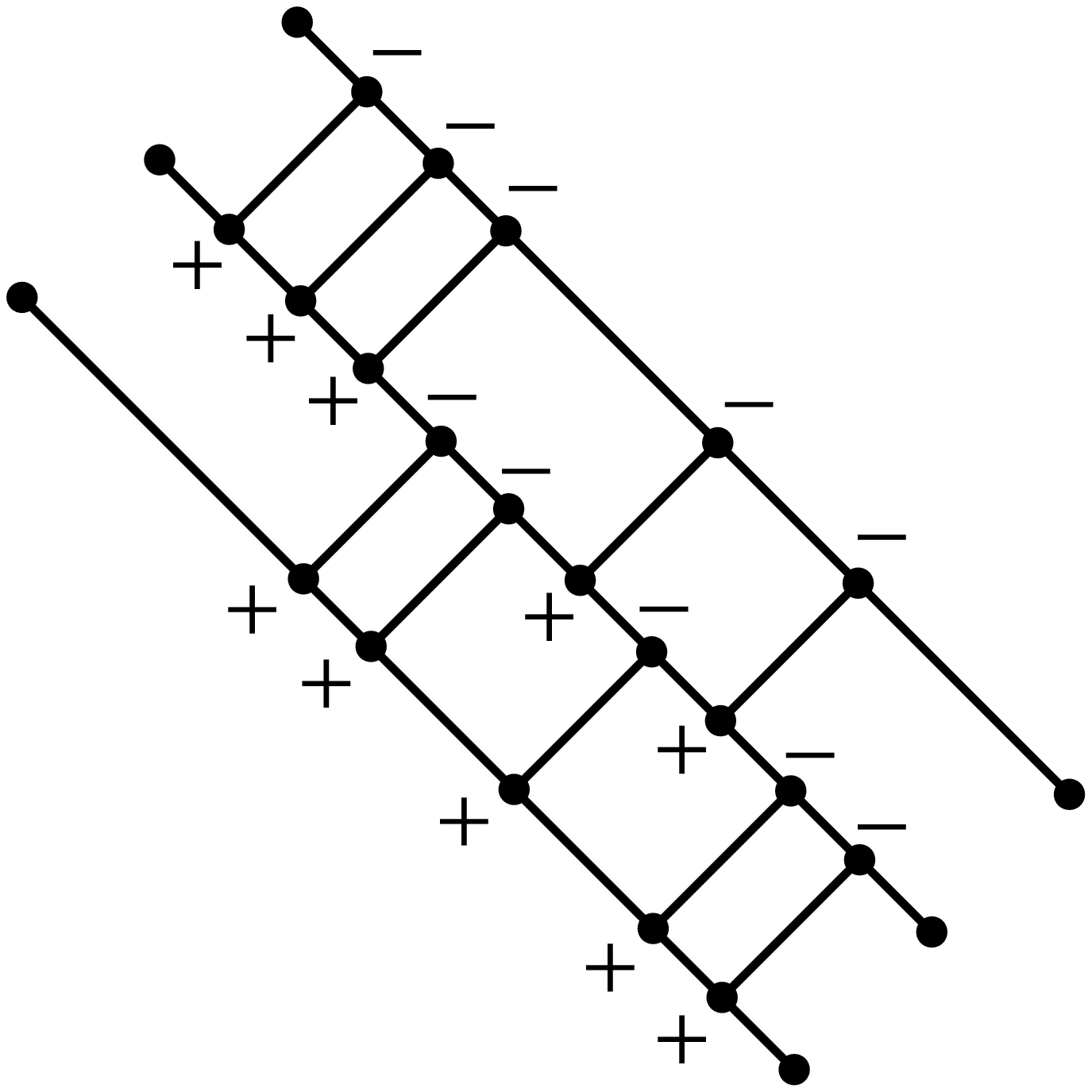}
\end{center}
\caption{A closed positive braid and a graph divide representing it}
\label{fig;ex_pb}
\end{figure}

\begin{figure}[ht!]
\begin{center}
\includegraphics[height=2.5cm]{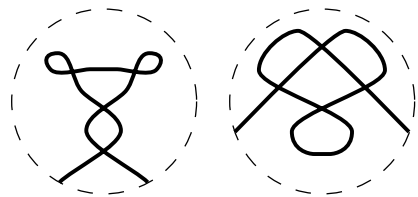}
\end{center}
\caption{Divides which represent the knot $10_{139}$}
\label{fig;10_139}
\end{figure}

\begin{remark}
In private communication, 
Ishikawa found that 
any closed positive braid can be represented 
as a graph divide link. 
For a given closed positive braid $L$, 
we draw a uni-trivalent graph 
which is homotopic to the canonical Seifert surface of $L$, 
and for endpoints of 
each edge corresponding to a positive crossing 
we give the different signs. 
For example, 
the knot $10_{139}$ in the table of Rolfsen \cite{rolfsen}
is the closure of the positive braid 
${\sigma _1}^2\sigma _2\sigma _1\sigma _2
{\sigma _1}^2{\sigma _2}^3$, 
where $\sigma _1$ and $\sigma _2$ are 
the canonical generators of the $3$--braid group. 
By this word, we draw a graph and give signs of vertices 
as illustrated in Figure \ref{fig;ex_pb}. 
The link of this graph divide is the knot $10_{139}$. 
Actually, the knot $10_{139}$ is the knot of each divide 
in Figure \ref{fig;10_139} 
\cite{acampo3, gibishi}. 
\end{remark}

\section{Four-dimensional invariants}\label{sect;d4inv}
%Slice euler characteristics of graph divide links.
%If a given graph is tree, 
%4-dimensional clasp number is determined. 

In this section, 
we determine some link invariants 
for graph divide links. 

Let $L$ be an oriented link, and 
let $F\subset B^4$ be a smooth, oriented 2--manifold 
with $\partial F=L$, 
where $B^4$ is the 4--ball bounded by $S^3$. 
We suppose that $F$ has no closed components, 
but $F$ is not assumed to be connected. 
We denote by $\chi _s(L)$ 
the greatest value of the euler characteristic $\chi (F)$
for such 2--manifolds $F\subset B^4$, 
and we call this invariant 
the \emph{slice euler characteristic}. 
In \cite{ru93}, Rudolph showed the following equality. 

\begin{theorem}{\rm\cite{ru93}}\label{thm;ru93}\qua
Let $L$ be 
the closure of a quasipositive braid 
with $n$ strings and $k$ quasipositive bands. 
Then we have $\chi _s(L)=n-k$. 
\qed
\end{theorem}

Combining Theorem \ref{thm;ru93} and 
the proof of Theorem \ref{thm;grdivqp}, 
we obtain the following result. 

\begin{proposition}\label{prop;slicegd}
Let $P=(G,\varphi )$ be a graph divide. 
Let $\delta (P)$ be the number of double points of $P$. 
Then the slice euler characteristic of the link of $P$ is 
$\chi (G)-2\delta (P)$. 
\end{proposition}

\begin{proof}
Let $P=(G,\varphi )$ be a graph divide 
and $\delta$ the number of double points of $P$. 
We may assume that 
$F(P;\{ \epsilon _x\} _{x\in V_P})$
in the proof of Theorem \ref{thm;grdivqp} 
is the projection image of a 2--manifold in 
$B^4-\{ \ast \}$ to $S^3$. 
We also denote this 2--manifold 
by $F(P;\{ \epsilon _x\} _{x\in V_P})$. 
By Theorem \ref{thm;ru93} 
and the construction of 
$F(P;\{ \epsilon _x\} _{x\in V_P})$ 
we have 
\[
\chi _s(
L(P;\{ \epsilon _x\} _{x\in V_P})
)=
\chi (
F(P;\{ \epsilon _x\} _{x\in V_P})
). \]
By means of constructions of an immersed 2--manifold 
$F(P;\{ \epsilon _x\} _{x\in V_P})$,  
we have 
\[
\chi (
F(P;\{ \epsilon _x\} _{x\in V_P})
)
=\chi (G)-2\delta (P). \]
\end{proof}

\begin{remark}
The slice euler characteristic of 
the link of graph divide 
does not depend on signs of vertices 
by the above proposition. 
\end{remark}

Let $L$ be an $r$--component link. 
The \emph{four-dimensional clasp number} of $L$ is 
the minimum number of the double points 
for transversely immersed $r$ disks in $B^4$ 
with boundary $L$ and with only finite double points 
as singularities \cite{ka3}. 
The author showed in \cite{ka3} that 
the inequality $c_s(L)\geq (r-\chi _s(L))/2$ holds 
for any link $L$.
In particular, 
if $G$ is a disjoint union of intervals and trees, 
the 4--dimensional clasp number of 
the link of the graph divide $P=(G,\varphi )$ 
is determined as below. 

\begin{corollary}\label{cor;d4clasp}
Let $P=(G,\varphi )$ be a graph divide. 
Let $\delta (P)$ be the number of double points of $P$. 
If $G$ is a disjoint union of intervals and trees, 
then the four-dimensional clasp number of the links of $P$ is 
$\delta (P)$. 
\end{corollary}

\begin{figure}[ht!]
\begin{center}
\includegraphics[height=2cm]{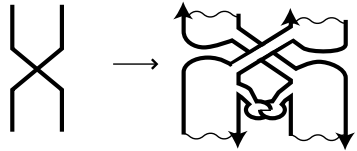}
\end{center}
\caption{Construction of immersed disks}
\label{fig;clasp}
\end{figure}

\begin{figure}[ht!]
\begin{center}
\includegraphics[height=5cm]{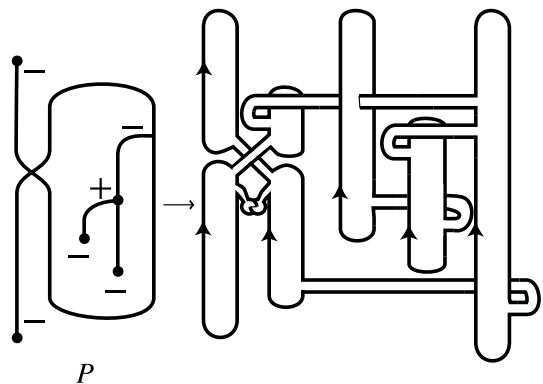}
\end{center}
\caption{An example of the immersed disk}
\label{fig;exclasp}
\end{figure}

\begin{proof}
By the above, 
the 4--dimensional clasp number of the links of $P$ is 
not less than $\delta (P)$. 
We show that 
the 4--dimensional clasp number of the links of $P$ is 
not greater than $\delta (P)$. 

For a given tree divide $P$ 
and signs of vertices $\{ \epsilon _x\} _{x\in V_P}$, 
we construct immersed disks in the 3--sphere 
by the almost same argument 
as that for the 2--manifold $F(P; \{ \epsilon _x\} _{x\in V_P})$ 
in Section \ref{sect;qpos}. 
For $P$ except near double points, 
we construct same parts of $F(P; \{ \epsilon _x\} _{x\in V_P})$. 
Around each double point of $P$, 
we construct parts of disks with a clasp 
as illustrated in Figure \ref{fig;clasp}. 
Non-vertical narrow bands may intersect some disks 
as explained in Step 3, 4, and 5 
of the construction of $F(P;\{ \epsilon _x\} _{x\in V_P})$. 
As an example the immersed disk in Figure \ref{fig;exclasp} 
is obtained from the pair of $P$ and $F(P;\{ \epsilon _x\} _{x\in V_P})$ 
in Figure \ref{fig;exFP} by the above algorithm. 
We regard the obtained diagram of immersed disks 
as a diagram of disks immersed in $B^4$ 
with boundary $L(P; \{ \epsilon _x\} _{x\in V_P})$. 
Then we have 
$c_s(L(P; \{ \epsilon _x\} _{x\in V_P}))\leq \delta (P)$. 
\end{proof}

\begin{remark}
In Example \ref{ex;fish2}, we represent the knot $8_{21}$ 
as a graph divide knot. 
It is known that 
the genus of this knot is $2$. 
In \cite{shi} Shibuya showed that 
for any link, 
the 3--dimensional clasp number is not less than 
the genus. 
Then the 3--dimensional clasp number of $8_{21}$ 
is not less than 2, and actually equal to $2$.  
By Corollary \ref{cor;d4clasp}, 
we have $c_s(8_{21})=1$. 
Therefore Corollary \ref{cor;d4clasp} cannot determine 
the 3--dimensional clasp numbers of graph divide links 
generally. 
\end{remark}

\begin{remark}
A'Campo \cite{acampo}, Gibson and Ishikawa \cite{gibishi} showed 
that 
if $P$ is a divide or a free divide without circle branches, 
the unknotting number of the link of $P$ is also equal to 
the number of double points of $P$.  
In \cite{gibishi}, Gibson and Ishikawa introduced 
a good unknotting operation for $P$. 
Furthermore in \cite{gib;tree} 
Gibson showed that 
the unknotting number of the link of a tree divide 
with some conditions 
is equal to the number of double points. 
\end{remark}

In \cite{gibishi},  
Gibson and Ishikawa checked that
there exist free divide links with no divide representation 
and that
the knot of the free divide with only one double point 
must be the trefoil or the mirror image of the knot $5_2$ 
in the table due to Rolfsen \cite{rolfsen}.  
Then Corollary \ref{cor;d4clasp} implies that 
each of the knots of Example \ref{ex;fish1} for $n\geq 0$ 
and Example \ref{ex;fish2} 
cannot be represented by any free divide (cf.\ \cite{gib;tree}). 
Therefore the class of graph divide links is 
a real extension of 
the class of free divide links. 

Then we are interested in the problem of 
the existence of the quasipositive link 
with no graph divide representation. 
By Proposition \ref{prop;slicegd}, 
we obtain the following result. 

\begin{theorem}\label{thm;qpnotgd}
If a graph divide knot $K$ is slice, 
then $K$ is trivial. 
Therefore, there exists a quasipositive link 
which is not a graph divide link. 
\end{theorem}

In the proof of the above theorem, 
we use the following lemma. 

\begin{lemma}\label{lem;embedtree}
If a graph $G$ is a connected tree and 
the map $\varphi \co G \rightarrow D$ is 
an embedding, 
then the link of $P=(G,\varphi )$ is a trivial knot. 
\end{lemma}

\begin{proof}
By moves illustrated in Figure \ref{fig;trans}, 
$G$ can be transformed to an arc $A$ embedded in $D$. 
By Lemma \ref{lem;trans}, 
the link of $P=(G,\varphi )$ is ambient isotopic to 
the link of $A$. 
\end{proof}

\begin{proof}[Proof of Theorem \ref{thm;qpnotgd}]
We suppose that 
$K=L((G,\varphi );\{ \epsilon _x\} _{x\in V_P})$ 
is a slice knot. 
The graph $G$ is connected since $K$ is a knot. 
Hence the euler characteristic of $G$ is not greater than $1$. 
By Proposition \ref{prop;slicegd}, 
the slice euler characteristic of $K$ is $\chi (G)-2\delta$ 
where $\delta$ is the number of double points of $P$. 
We have $\chi (G)=2\delta +1$ since the knot $K$ is slice. 
Then the euler characteristic of $G$ must be $1$ and 
$\delta$ must be $0$. 
Therefore $G$ is a connected tree and 
$\varphi$ is an embedding. 
By Lemma \ref{lem;embedtree}, $K$ is a trivial knot. 

The mirror image of the knot $8_{20}$ 
in the table due to Rolfsen \cite{rolfsen} is 
the closure of the quasipositive braid 
$\sigma _1^3\sigma _2\sigma _1^{-3}\sigma _2$, 
where $\sigma _1$ and $\sigma _2$ are canonical generators of 
the 3--braid group. 
This knot is non-trivial and slice. 
Then it is an example of a quasipositive link 
which is not a graph divide link. 
\end{proof}

\begin{remark}
A {\em strongly quasipositive braid} is 
the product of positive bands. 
A {\em positive band} is a braid 
$\sigma _{ij}
=\sigma _i\cdots \sigma _{j-1}\sigma _{j}
{\sigma _{j-1}}^{-1}\cdots {\sigma _i}^{-1}$,
where $\sigma _k$'s are canonical generators of the braid group 
and $i$ is less than $j$. 
A {\em strongly quasipositive link} is 
the closure of a strongly quasipositive braid. 
For any strongly quasipositive link, 
the slice euler characteristic is equal to 
the euler characteristic. 
Hirasawa recently showed that 
any divide link is strongly quasipositive \cite{hira;strqp}. 
The knot of Example \ref{ex;fish2} is fibered, 
but not strongly quasipositive, 
because 
the slice euler characteristic is $-1$ 
but the euler characteristic is $-3$, so 
the slice euler characteristic is not equal to 
the euler characteristic. 
\end{remark}

\Addresses\recd
\end{document}

%% file: agtout.tex
\def\ifplaintex{\expandafter\ifx\csname documentclass\endcsname\relax}

\def\gtp{{\mathsurround=0pt\it $\cal G\mskip-2mu$eometry \&\ 
$\cal T\!\!$opology $\cal P\!$ublications}}  % GT publications

\def\recd{{\small Received:\qua\receiveddate\ifx\reviseddate\relax
\else\qquad Revised:\qua\reviseddate\fi\par}} 

\def\lognumber#1{\def\thelognumber{#1}}
\def\volumenumber#1{\def\thevolumenumber{#1}}
\def\volumeyear#1{\def\thevolumeyear{#1}}
\def\papernumber#1{\def\thepapernumber{#1}}
\def\pagenumbers#1#2{\def\startpage{#1}\def\finishpage{#2}}
\def\published#1{\def\publishdate{#1}}

\def\received#1{\def\receiveddate{#1}}
\def\revised#1{\def\reviseddate{#1}}
\def\accepted#1{\def\accepteddate{#1}}

\let\\\par\let\thelognumber\relax\let\thevolumenumber\relax
\let\thepapernumber\relax\let\thevolumeyear\relax\let\startpage\relax
\let\finishpage\relax\let\publishdate\relax\let\receiveddate\relax
\let\reviseddate\relax\let\accepteddate\relax\let\theasciititle\relax
\let\theasciiauthors\relax
\let\theasciiabstract\relax

\let\theasciiemail\relax

\ifplaintex
\font\logobig=cmssbx10 scaled 3836
\font\logomed=cmssbx10 scaled 2557
\else
\font\logobig=cmssbx10 scaled 4200
\font\logomed=cmssbx10 scaled 2800
\fi

\long\def\makeagttitle{   %%% start of definition of \makeagttitle
\count0=\startpage
\agt\hfill      %   Journal title (top left) 
\hbox to 45truept{\vbox to 0pt{\vglue -13truept{\logomed A\kern -.37em{\logobig 
T}\kern -.38em G}\vss}\hss}
\break
{\small Volume \thevolumenumber\ (\thevolumeyear)
\startpage--\finishpage\nl
Published: \publishdate}

\vglue .25truein

{\parskip=0pt\leftskip 0pt plus
1fil\def\\{\par\smallskip}{\Large\bf\thetitle}\par\medskip} \vglue
0.05truein

{\parskip=0pt\leftskip 0pt plus 1fil\def\\{\par}{\sc\theauthors}
\par\medskip}%
 
\vglue 0.03truein 

{\small\leftskip 25truept\rightskip 25truept{\bf Abstract}\stdspace\theabstract

{\bf AMS Classification}\stdspace\theprimaryclass
\ifx\thesecondaryclass\relax\else; \thesecondaryclass\fi\par
{\bf Keywords}\stdspace \thekeywords\par}\vglue 7truept

}   %%%% end of definition of \makeagttitle

\ifplaintex
\font\phead=cmsl9 scaled 950
\font\pnum=cmbx10 scaled 913
\font\pfoot=cmsl9 scaled 950
\headline{\vbox to 0pt{\vskip -4.5mm\line{\small\phead\ifnum
\count0=\startpage ISSN 1472-2739 (on-line) 1472-2747 (printed)
\hfill {\pnum\folio}\else\ifodd\count0\def\\{ }% 
\ifx\theshorttitle\relax\thetitle\else\theshorttitle\fi\hfill{\pnum\folio}
\else\def\\{ and }{\pnum\folio}\hfill\ifx\theshortauthors\relax\theauthors
\else\theshortauthors\fi\fi\fi}\vss}}
\footline{\vbox to 0pt{\vglue 0mm\line{\small\pfoot\ifnum\count0=\startpage
\copyright\ \gtp\hfill\else
\agt, Volume \thevolumenumber\ (\thevolumeyear)\hfill\fi}\vss}}
\else
\font=cmsl9 scaled 1050
\font\lnum=cmbx10 
\font=cmsl9 scaled 1050
\makeatletter
\def\@oddhead{{\small\lhead\ifnum\count0=\startpage ISSN 1472-2739 
(on-line) 1472-2747 (printed)\hfill {\lnum\number\count0}\else\ifodd\count0
\def\\{ }\ifx\theshorttitle\relax \thetitle \else\theshorttitle\fi\hfill
{\lnum\number\count0}\else\def\\{ and }{\lnum\number\count0}
\hfill\ifx\theshortauthors\relax 
\theauthors\else\theshortauthors\fi\fi\fi}}\def\@evenhead{\@oddhead}
\def\@oddfoot{\small\lfoot\ifnum\count0=\startpage\copyright\ \gtp\hfill\else
\agt, Volume \thevolumenumber\ (\thevolumeyear)\hfill\fi}
\def\@evenfoot{\@oddfoot}
\makeatother
\fi
\let\maketitlepage\makeagttitle
\let\makeshorttitle\maketitlepage
\let\maketitle\maketitlepage

\newwrite\gtoutfile
\long\gdef\makeheadfile{  %%% start of definition of \makeheadfile
{\def\\{, }\def\s{ }
\immediate\openout\gtoutfile head.xxx
\immediate\write\gtoutfile{Proxy-for: \ifx\theasciiauthors\relax
\theauthors\else\theasciiauthors\fi\s<\ifx\theasciiemail\relax\theemail\else\theasciiemail\fi>}
\immediate\write\gtoutfile{\noexpand\\}
\immediate\write\gtoutfile{Authors: \ifx\theasciiauthors\relax
\theauthors\else\theasciiauthors\fi}
{\def\\{ }\immediate\write\gtoutfile{Title: \ifx\theasciititle\relax
\thetitle\else\theasciititle\fi}}
\immediate\write\gtoutfile{Subj-class: GT or SG, GR etc}
\immediate\write\gtoutfile{MSC-class: \theprimaryclass\ifx\thesecondaryclass\relax\else, \thesecondaryclass\fi}
\immediate\write\gtoutfile{Journal-ref: Algebr. Geom. Topol. \thevolumenumber\s
(\thevolumeyear) \startpage-\finishpage}
\immediate\write\gtoutfile{Comments: Published by Algebraic and
Geometric Topology at}
\immediate\write\gtoutfile{\s\s\s  http://www.maths.warwick.ac.uk/agt/AGTVol\thevolumenumber/agt-\thevolumenumber-\thepapernumber.abs.html}
\immediate\write\gtoutfile{\noexpand\\}
\immediate\write\gtoutfile{}
\ifx\theasciiabstract\relax
\immediate\write\gtoutfile{\theabstract}\else
\immediate\write\gtoutfile{\theasciiabstract}\fi
\immediate\write\gtoutfile{}
\immediate\write\gtoutfile{\noexpand\\}
\immediate\write\gtoutfile{}
\immediate\closeout\gtoutfile}}  %%% end of definition of \makeheadfile

\def\maketitlepage{\makeagttitle\makeheadfile}
\let\makeshorttitle\maketitlepage
\let\maketitle\maketitlepage

%% file: agt-4-26.bbl
\begin{thebibliography}


 \bibitem{acampo1}
     {\bf N A'Campo}, 
     {\em Real deformations and complex topology 
      of plane curve singularities}, 
     Annales de la Facult\'{e} des Sciences de Toulouse
     8 (1999) 5--23    \MR{1721511}
 \bibitem{acampo}
     {\bf N A'Campo}, 
     \textit{Generic immersions of curves, knots, 
      monodromy and gordian number}, 
     Publ. Math. IHES 88 (1998) 151--169 \MR{1733329}
 \bibitem{acampo3}
     {\bf N A'Campo}, 
     \textit{Planer trees, slalom curves and hyperbolic knots}, 
     Publ. Math. IHES 88 (1998) 171--180 \MR{1733330} 
 \bibitem{BO}
      {\bf M Boileau}, {\bf S Orevkov}, 
      {\em Quasipositivit\'{e} d'une courbe analytique 
           dans une boule pseudo-convexe}, 
      C. R. Acad. Sci. Paris 332 (2001) 825--830 \MR{1836094} 
 \bibitem{gib;tree}
     {\bf W Gibson}, 
     {\em Links and Gordian numbers associated with 
          generic immersions of trees}, 
     Proceedings of Art of Low Dimensional Topology V$\!$I$\!$I (2001) 
     27--35 
 \bibitem{gibishiOD}
     {\bf W Gibson}, {\bf M Ishikawa}, 
     {\em Links of oriented divides 
          and fibrations in link exteriors}, 
     Osaka J. Math. 39 (2002) 681--703 \MR{1932288}
 \bibitem{gibishi}
     {\bf W Gibson}, {\bf M Ishikawa}, 
     {\em Links and Gordian numbers associated with 
          generic immersions of intervals}, 
     Topol. Appl. 123 (2002) 609--636 \MR{1924054}
 \bibitem{hira}
     {\bf M Hirasawa}, 
     {\em Visualization of A'Campo's fibered links 
             and unknotting operations}, 
     Topol. Appl. 121 (2002) 287--304 \MR{1903697}
 \bibitem{hira;strqp}
     {\bf M Hirasawa}, 
     in preparation 
 \bibitem{ishi;Jpoly}
     {\bf M Ishikawa}, 
     {\em The $\mathbb{Z}_2$--coefficient 
          Kauffman state model on divides}, 
     preprint 
 \bibitem{ka3}
     {\bf T Kawamura}, 
     {\em On unknotting numbers and 
      four-dimensional clasp numbers of links}, 
     Proc. Amer. Math. Soc. 130 (2002) 243--252 \MR{1855642} 
 \bibitem{ka6}
     {\bf T Kawamura}, 
     {\em Quasipositivity of links of divides and free divides},  
     Topol. Appl. 125 (2002) 111--123 \MR{1931179}
 \bibitem{rolfsen}
     {\bf D Rolfsen}, {\em Knots and links}, 
     Mathematics Lecture Series  7, 
     Publish or Perish, Berkeley, Calif. (1976) \MR{0515288} 
  \bibitem{rualg}
     {\bf L Rudolph}, 
     {\em Algebraic functions and closed braids}, 
     Topology 22 (1983) 191--201 \MR{0683760} 
 \bibitem{rutop90}
     {\bf L Rudolph}, 
     {\em Totally tangential links of intersection 
      of complex plane curves with round spheres}, 
     from: ``Topology '90 (Columbus OH (1990)'', 
     de Gruyter, Berlin (1992) 343--349 \MR{1184419}
 \bibitem{ru93}
     {\bf L Rudolph}, 
     {\em Quasipositivity as an obstruction to sliceness}, 
     Bull. Amer. Math. Soc. 29 (1993) 51--59 \MR{1193540}
 \bibitem{shi}
     {\bf T Shibuya}, 
     {\em Some relations among various numerical invariants 
          for links}, 
     Osaka J. Math. 11 (1974) 313--322 \MR{0353295}

\end{thebibliography}
